\documentclass[11pt,a4paper,leqno]{amsart}

\usepackage{amsmath,amssymb, amsfonts}

\def\<#1,#2>{\langle\,#1,\,#2\,\rangle}

\def\qed{\ensuremath{\hfill\Box}}
\def\qed{\ensuremath{\hfill\Box}}

\def\f{\varphi}
\def\phi{\varphi}

\theoremstyle{plain}
\newtheorem{Lemma}{Lemma}[section]
\newtheorem{Proposition}[Lemma]{Proposition}
\newtheorem{Theorem}[Lemma]{Theorem}

\theoremstyle{definition}
\newtheorem{Definition}[Lemma]{Definition}

\theoremstyle{remark}
\newtheorem{Remark}[Lemma]{Remark}

\title{Remarks on the product of harmonic forms}
\author{Liviu Ornea}
\thanks{Both authors are partially supported by CNCSIS grant PNII IDEI contract 529/2009. The second-named author acknowledges also partial support from SFB/TR 12.}

\address{University of Bucharest, Faculty of Mathematics, 14
Academiei str., 70109 Bucharest, Romania \emph{and} 
Institute of Mathematics ``Simion Stoilow" of the Romanian Academy, 
21, Calea Grivitei str.
010702-Bucharest, Romania}
\email{lornea@gta.math.unibuc.ro, liviu.ornea@imar.ro}

\author{Mihaela Pilca}
\address{Mathematisches Institut, Universit\"at zu K\"oln,
Weyertal 86-90 D-50931 K\"oln, Germany \emph{and} 
Institute of Mathematics ``Simion Stoilow" of the Romanian Academy, 
21, Calea Grivitei str.
010702-Bucharest, Romania}
\email{mpilca@mi.uni-koeln.de}

\begin{document}

\begin{abstract}

A metric is formal if all products of harmonic forms are again harmonic. The existence of a formal metric implies Sullivan formality of the manifold, and hence formal metrics can exist only in presence of a very restricted topology. We show that a warped product metric is formal if and only if the warping function is constant and derive further topological obstructions to the existence of formal metrics. In particular, we determine necessary and sufficicient conditions for a Vaisman metric to be formal.

\medskip


\smallskip

\end{abstract}

\maketitle

\section{Introduction}

A fundamental problem in algebraic topology is the reading of the homotopy type of a space in terms of cohomological data. A precise definition of this property was given by Sullivan in \cite{su} and called {\em formality}. As concerns manifolds, it is known {\em e.g.} that all compact Riemannian symmetric spaces and all compact K\"ahler manifolds are formal. For a recent survey of topological formality, see \cite{ps}.

Sullivan also observed that if a compact manifold admits a metric such that the wedge product of any two harmonic forms is again harmonic, then, by Hodge theory,  the manifold is formal. This motivated Kotschick to give the following:

\begin{Definition}(\cite{kot})
A Riemannian metric is called {\em (metrically) formal} if all
wedge products of harmonic forms are harmonic.

A closed manifold is called {\em geometrically formal} if it admits a formal
Riemannian metric.
\end{Definition}

In particular, the length of any harmonic form with respect to a formal metric is (pointwise) constant. This larger class of metrics having all harmonic (1-)forms of constant length naturally appears in other geometric contexts, for instance in the study of certain systolic inequalities, and has been investigated in \cite{nagy}, \cite{nagyver}.

Classical examples of geometrically formal manifolds are compact symmetric spaces. In \cite{kot_terz1} and \cite{kot_terz2} more general examples are provided, both of geometrically formal and of formal but non-geometrically formal homogeneous manifolds.

Geometric formality imposes strong restrictions on the (real) cohomology of the manifold. For example,  it is proven in \cite{kot} that a manifold admits a non-formal metric if and only if it is not a rational homology sphere.

In this note, we shall obtain further obstructions to formality. We shall see that if a compact manifold with $b_1=p\geq 1$ admits a formal metric, and if there exist two vanishing Betti numbers such that the distance between them is not larger than $p+2$, then all the intermediary Betti numbers must be zero too. Also, a conformal class of metrics on an even-dimensional compact manifold with non-zero middle Betti number can contain no more than one formal metric. 

Our main concern will be the formality of warped products. We shall show that a warped product metric on a compact manifold is formal if and only if the warping function is constant. On the way, we shall also provide a proof for the known fact (stated for instance in \cite{kot}) that a product of formal metrics is formal. 

Unlike K\"ahler manifolds, which are known to be formal, for the time being, nothing is known about the Sullivan formality of locally conformally K\"ahler  (in particular Vaisman) manifolds. In the last section of this note, we shall discuss compact Vaisman manifolds, whose universal cover is a special type of warped product, a Riemannian cone to be precise, and we shall find obstructions to the metric formality of a Vaisman metric. 

\smallskip

We end this introduction with the following straightforward, but useful characteri\-sation of geometric formality:
\begin{Lemma}\label{lharm}
Let $\alpha$ and $\beta$ be two harmonic forms on a compact Riemannian manifold $(M^n,g)$. Then $\alpha\wedge\beta$ is harmonic if and only if the following equality is satisfied:
\begin{equation}\label{cdtharm}
\sum_{i=1}^{n}(e_i\lrcorner \alpha)\wedge\nabla_{e_i}\beta=-(-1)^{|\alpha||\beta|}\sum_{i=1}^{n}(e_i\lrcorner \beta)\wedge\nabla_{e_i}\alpha,
\end{equation}
where $\{e_i\}_{i=\overline{1,n}}$ is a local orthonormal basis of vector fields. Thus, the metric $g$ is formal if and only if \eqref{cdtharm} holds  for any two $g$-harmonic forms.
\end{Lemma}
\begin{proof}
Since we are on a compact manifold, a differential form is harmonic if and only if it is closed and coclosed. As $\alpha\wedge\beta$ is closed, we have to show that \eqref{cdtharm} is equivalent to $\alpha\wedge\beta$ being coclosed. This is implied by the following:
\begin{equation*}
\begin{split}
&\delta(\alpha\wedge\beta)=-\sum_{i=1}^{n}e_i\lrcorner \nabla_{e_i}(\alpha\wedge\beta)=-\sum_{i=1}^{n}e_i\lrcorner (\nabla_{e_i}\alpha\wedge\beta+\alpha\wedge\nabla_{e_i}\beta)\\
&=\delta\alpha\wedge\beta-(-1)^{|\alpha|}\sum_{i=1}^{n}\nabla_{e_i}\alpha\wedge(e_i\lrcorner\beta)-\sum_{i=1}^{n}(e_i\lrcorner\alpha)\wedge\nabla_{e_i}\beta+(-1)^{|\alpha|}\alpha\wedge\delta\beta\\
&=-(-1)^{|\alpha||\beta|}\sum_{i=1}^{n}(e_i\lrcorner\beta)\wedge\nabla_{e_i}\alpha-\sum_{i=1}^{n}(e_i\lrcorner\alpha)\wedge\nabla_{e_i}\beta.
\end{split}
\end{equation*}

$\qed$ 

\end{proof}

\section{Geometric formality of warped product metrics}

\subsection{Riemannian products}

For the sake of completeness and as a first step in the study of geometrically formal warped products, we provide a proof for the formality of Riemannian product metrics.

Let $(M^{n+m},g)=(M^n_1,g_1)\times (M^m_2,g_2)$ be the Riemannian product of two compact manifolds and let $\pi_i:M\to M_i$ be the natural projections, which are totally geodesic Riemannian submersions.

One may describe the bundle of $p$-forms on $M$ as follows:
\begin{equation}\label{descompformp}
\Lambda^p M=\overset{p}{\underset{k=0}{\bigoplus}}\pi_1^*(\Lambda^k M_1)\otimes \pi_2^*(\Lambda^{p-k} M_2).
\end{equation}
This identification also works for the space of harmonic forms, namely the harmonic forms  on $(M,g)$ can be described in terms of the harmonic forms on the factors $(M_1,g_1)$ and $(M_2,g_2)$. To this end let  $\mathcal{H}^k(M_i,g_i)$ be the space of harmonic $k$-forms on $M_i$ and let $b_k(M_i)$ be the Betti numbers of $M_i$, $i=1,2$.  We can prove:
\begin{Lemma}\label{lbas}
Let  $\{\alpha^k_{1}, \dots, \alpha^k_{b_k(M_1)}\}$ (resp. $\{\beta^k_{1}, \dots, \beta^k_{b_k(M_2)}\}$) be a basis of $\mathcal{H}^k(M_1,g_1)$ (resp. $\mathcal{H}^k(M_2,g_2)$). Then the forms:
\begin{equation}\label{basisharm}
\{\pi_1^*(\alpha^k_s)\wedge\pi_2^*(\beta^{p-k}_l) |\; 1\leq s\leq b_k(M_1), 1\leq l\leq b_{p-k}(M_2), 0\leq k\leq p\}
\end{equation}
form a basis of the space of  $\mathcal{H}^{p}(M,g)$, for each $0\leq p\leq m+n$.
\end{Lemma}

\begin{proof} By Hodge theory for compact manifolds, the dimension of the space of harmonic $p$-forms is equal to the $p^{\text{th}}$ Betti number, and hence:
\begin{equation*}
\begin{split}
\dim(\mathcal{H}^{p}(M,g))&=b_p(M)=\sum_{k=0}^{p}b_{k}(M_1)b_{p-k}(M_2)\\
&=\sum_{k=0}^{p}\dim(\mathcal{H}^{k}(M_1,g_1))\dim(\mathcal{H}^{p-k}(M_2,g_2)).
\end{split}
\end{equation*}
It sufficies to show that the forms in \eqref{basisharm} are $g$-harmonic (since they are linear independent and in the right number, they build a basis of $\mathcal{H}^{p}(M,g)$). It is enough to check that each form $\pi_1^*(\alpha)\wedge\pi_2^*(\beta)$ is $g$-harmonic if $\alpha$ is a $g_1$-harmonic form on $M_1$ and $\beta$ is a $g_2$-harmonic form on $M_2$. We first show that $\pi_1^*(\alpha)$ and $\pi_2^*(\beta)$ are $g$-harmonic forms on $M$, then use Lemma \ref{lharm}. Since the manifolds are compact, a form is harmonic if and only if it is closed and coclosed.

As the exterior differential commutes with the pull-back,  $\pi_1^*(\alpha)$ and $\pi_2^*(\beta)$ are closed forms on $M$.

For the codifferential $\delta$ on $M$ we first check that it commutes with the pull-back given by the projections $\pi_i$. Let $\{e_i\}_{i=\overline{1,n}}$ be a local orthonormal basis on $M_1$ and $\{f_j\}_{j=\overline{1,m}}$ be a local orthonormal basis on $M_2$, which we lift to $M$ and thus obtain a local orthonormal basis of $M$ : $\{\tilde{e_i},\tilde{f_j}\}_{i=\overline{1,n}; j=\overline{1,m}}$. We consider the following decomposition of the codifferential on $M$: $\delta=\delta_1+\delta_2$, where
\[\delta_1:=-\sum_{i=1}^{n}{\tilde{e_i}}\lrcorner\nabla_{\tilde{e_i}}, \quad \delta_2:=-\sum_{j=1}^{m}{\tilde{f_j}}\lrcorner\nabla_{\tilde{f_j}}.\]
The pull-back of any $p$-form $\alpha$ on $M_1$ is automatically in the kernel of $\delta_2$  since $\nabla_{X}(\pi_1^*(\alpha))=0$, for any vector field $X$ tangent to $M_2$:
\begin{equation*}
\begin{split}
(\nabla_{X}(\pi_1^*(\alpha)))(&Y_1,\dots,Y_p)=X(\alpha({\pi_1}_*(Y_1), \dots, {\pi_1}_*(Y_p))\circ\pi_1)\\
&-\sum_{j=1}^{p} \alpha({\pi_1}_*(Y_1),\dots, {\pi_1}_*(\nabla_X Y_j),\dots,{\pi_1}_*(Y_p))\circ\pi_1=0,
\end{split}
\end{equation*}
where $\{Y_i\}_{i=\overline{1,p}}$ are any tangent vector fields to $M$.

We then obtain:
\begin{equation}\label{codiffa}
\begin{split}
\delta(\pi_1^*(\alpha))&=\delta_1(\pi^*_1(\alpha))=-\sum_{i=1}^{n}{\tilde{e_i}}\lrcorner\nabla_{\tilde{e_i}}(\pi_1^*(\alpha))\\
&=-\sum_{i=1}^{n}{\tilde{e_i}}\lrcorner\pi_1^*(\nabla^{g_1}_{e_i}\alpha)=-\sum_{i=1}^{n}\pi_1^*(e_i\lrcorner\nabla^{g_1}_{e_i}\alpha)=\pi_1^*(\delta^{g_1}\alpha),
\end{split}
\end{equation}
where by $\nabla^{g_1}$ we denote the Levi-Civita connection and by $\delta^{g_1}$ the codifferential of $g_1$ on $M_1$. Since the roles of $M_1$ and $M_2$ are symmetric, we obtain a similar commutation relation for any form $\beta$ on $M_2$:
\begin{equation}\label{codiffb}
\delta(\pi_2^*(\beta))=\pi_2^*(\delta^{g_2}\beta).
\end{equation}

From \eqref{codiffa}, \eqref{codiffb} and the closedness of $\pi_1^*(\alpha)$ and $\pi_2^*(\beta)$, it follows that $\pi_1^*(\alpha)$ and $\pi_2^*(\beta)$ are $g$-harmonic, if $\alpha$ and $\beta$ are $g_1$- respectively $g_2$-harmonic forms. 

By Lemma \ref{lharm}, in order to show that $\pi_1^*(\alpha)\wedge\pi_2^*(\beta)$ is harmonic, we have to check that condition \eqref{cdtharm} is fulfilled. Considering again an adapted local othonormal basis $\{\tilde{e_i},\tilde{f_j}\}_{i=\overline{1,n}; j=\overline{1,m}}$ as above, it follows that \eqref{cdtharm} holds for $\pi_1^*(\alpha)$ and $\pi_2^*(\beta)$, since in each term one factor vanishes: $\tilde{e_i}\lrcorner\pi_1^*(\beta)=0$ and $\nabla_{\tilde{e_i}}(\pi_1^*(\beta))=0$, for all $i=1,\dots, n$ and $\tilde{f_j}\lrcorner\pi_2^*(\alpha)=0$ and $\nabla_{\tilde{f_j}}(\pi_2^*(\alpha))=0$, for all $j=1,\dots,m$. 

$\qed$ 
\end{proof}

We are now ready to prove:

\begin{Proposition}\label{prod}
If $(M_1,g_1)$ and $(M_2,g_2)$ are two compact \mbox{Riemannian} manifolds with formal metrics, then the metric $g=g_1+g_2$ on the product manifold $M=M_1\times M_2$ is also formal.
\end{Proposition}

\begin{proof}
Let $\gamma\in\Omega^p{M}$ and $\gamma'\in\Omega^q{M}$ be two harmonic forms on $M$. By Lemma \ref{lbas}, $\gamma$ and $\gamma'$ are given by linear combinations with real coefficients of the basis elements in \eqref{basisharm}. Thus, it is enough to check that the exterior product of any two such basis elements is a harmonic form on $M$. But:
\begin{equation*}
\begin{split}
(\pi_1^*(\alpha)\wedge\pi_2^*(\beta))\wedge (\pi_1^*(\alpha')\wedge\pi_2^*(\beta'))=(-1)^{|\alpha'||\beta|}\pi_1^*(\alpha\wedge\alpha')\wedge\pi_2^*(\beta\wedge\beta'),
\end{split}
\end{equation*}
which is $g$-harmonic on $M$ by Lemma \ref{lbas} and by the formality of $g_1$ and $g_2$ (as $\alpha\wedge\alpha'$ is again a $g_1$-harmonic form and $\beta\wedge\beta'$ a $g_2$-harmonic form). 

$\qed$ 
\end{proof}

\begin{Proposition}\label{uconf}
Let $M^{2n}$ be an even-dimensional compact manifold whose middle Betti number $b_n(M)$ is non-zero. Then, in any conformal class of metrics there is at most one formal metric (up to homothety).
\end{Proposition}

\begin{proof}
Let $[g]$ be a class of conformal metrics on $M$ and suppose there are two formal metrics $g_1$ and $g_2=e^{2f}g_1$ in $[g]$. The main observation is that in the middle dimension the kernel of the codifferential is invariant at conformal changes of the metric, so that there are the same harmonic forms for all metrics in a conformal class: $\mathcal{H}^n(M,g_1)=\mathcal{H}^n(M,g_2)$. As $b_n(M)\geq 1$ there exists a non-trivial $g_1$-harmonic (and thus also $g_2$-harmonic) $n$-form $\alpha$ on $M$. The length of $\alpha$ must then be constant with respect to both metrics, which are assumed to be formal and thus we get:
\[g_2(\alpha,\alpha)=e^{2nf}g_1(\alpha,\alpha),\]
which shows that $f$ must be constant. 

$\qed$ 
\end{proof}

\smallskip

Using the product construction to assure that the middle Betti number is non-zero, one can build such examples of formal metrics which are unique in their conformal class. 

Other examples are provided by manifolds with ``big'' first Betti number, as follows from the following property of ``propagation" of Betti numbers on geometrically formal manifolds proven in \cite[Theorem 7]{kot}:
{\em if $b_1(M)= p\geq 1$, then $b_q(M)\geq \binom{p}{q}$, for all $1\leq q\leq p$.} In particular, if $b_1(M^{2n})\geq n$, then $b_n(M^{2n})\geq 1$.

Another property of the Betti numbers of geometrically formal manifolds is given by:

\begin{Proposition}\label{bettzero}
Let $M^n$ be a compact geometrically formal manifold with $b_1(M)=p\geq 1$. 
If there exist two Betti numbers that vanish: $b_k(M)=b_{k+l}(M)=0$, for some $k$ and $l$ with $0<k+l<n$ and $0<l\leq p+1$, then all intermediary Betti numbers must vanish: $b_i(M)=0$, for $k\leq i\leq k+l$.
In particular, if there exists $k\geq \frac{n-p-1}{2}$ such that $b_k(M)=0$, then $b_i(M)=0$ for all $k\leq i\leq n-k$.
\end{Proposition}

\begin{proof}
Let $\{\theta_1,\dots,\theta_p\}$ be an orthogonal basis of $g$-harmonic $1$-forms, where $g$ is a formal metric on $M$. We first notice that here is no ambiguity in considering the orthogonality with respect to the global scalar product or to the pointwise inner product, because, when restricting ourselves to the space of harmonic forms of a formal metric, these notions coincide. This is mainly due to \cite[Lemma~4]{kot}, which states that the inner product of any two harmonic forms is a constant function. Thus, if two harmonic forms $\alpha$ and $\beta$ are orthogonal with respect to the global product, we get: $0=(\alpha,\beta)=\int_M<\alpha,\beta>dvol_g=<\alpha,\beta>vol(M)$, showing that their pointwise inner product is the zero-function.

It is enough to show that $b_{k+1}(M)=0$ and then use induction on $i$. Let $\alpha$ be a harmonic $(k+1)$-form. By formality, $\theta_1\wedge\theta_2\wedge\cdots\wedge\theta_{l-1}\wedge\alpha$ is a harmonic $(k+l)$-form and thus must vanish, since $b_{k+l}(M)=0$. On the other hand, $\theta_j^{\sharp}\lrcorner\alpha=(-1)^{k(n-k-1)}*(\theta_j\wedge*\alpha)$ is a harmonic $k$-form, again by formality. As $b_k(M)=0$, it follows that $\theta_j^{\sharp}\lrcorner\alpha=0$, for $1\leq j\leq p$. Then, using that $\{\theta_1,\dots,\theta_p\}$ are also orthogonal, we obtain:
\[0=\theta_1^{\sharp}\lrcorner \cdots \lrcorner \theta_{l-1}^{\sharp}\lrcorner(\theta_1\wedge\cdots\wedge\theta_{l-1}\wedge\alpha)=\pm|\theta_1|^2\cdots|\theta_{l-1}|^2\alpha,\]
which implies that $\alpha=0$, because each $\theta_j$ has non-zero constant length. This shows that $b_{k+1}(M)=0$.

$\qed$ 
\end{proof}

\subsection{Warped products}

We now pass to the setting we are mainly interested in, namely the warped products. 

Let $(B^n,g_B)$ and $(F^m,g_F)$ be two Riemannian manifolds and let $\f>0$ be a smooth function on $B$. Let $M=B\times_{\phi} F$ be the warped product with the metric $g=\pi^*(g_B)+(\phi\circ\pi)^2\sigma^*(g_F)$, where $\pi:M\to B$ and $\sigma:M\to F$ are  the natural projections. 

Let $\{e_i\}_{i=\overline{1,n}}$ be a local orthonormal basis on $B$ and let $\{f_j\}_{j=\overline{1,m}}$ be a local orthonormal basis on $F$, which we lift to $M$ and thus obtain a local orthonormal basis of $M$: $\displaystyle\{\tilde{e_i},\frac{1}{\phi\circ\pi}\tilde{f_j}\}_{i=\overline{1,n}; j=\overline{1,m}}$. 

Consider the following decomposition of the codifferential on $M$: $\delta=\delta_1+\delta_2$, where
\[\delta_1:=-\sum_{i=1}^{n}{\tilde{e_i}}\lrcorner\nabla_{\tilde{e_i}}, \quad \delta_2:=-\frac{1}{(\phi\circ\pi)^2}\sum_{j=1}^{m}{\tilde{f_j}}\lrcorner\nabla_{\tilde{f_j}}.\]

We first determine the commutation relations between the pull-back of forms on $B$ and $F$ with $\delta_1$ and $\delta_2$.
\begin{Lemma}
For $\alpha\in\Omega^{*}(B)$ and $\beta\in\Omega^{*}(F)$, the following relations hold:
\begin{equation}\label{del2pi}
\delta_2(\pi^*(\alpha))=-\frac{m}{\phi\circ\pi}\mathrm{grad}(\phi\circ\pi)\lrcorner\pi^*(\alpha),
\end{equation}
\begin{equation}\label{del1pi}
\delta_1(\pi^*(\alpha))=\pi^*(\delta^{g_B}(\alpha)),
\end{equation}
\begin{equation}\label{del1sig}
\delta_1(\sigma^*(\beta))=0,
\end{equation}
\begin{equation}\label{del2sig}
\delta_2(\sigma^*(\beta))=\frac{1}{(\phi\circ\pi)^2}\sigma^*(\delta^{g_F}(\beta)).
\end{equation}
\end{Lemma}

\begin{proof}
Let $\alpha\in\Omega^{p+1}(B)$. For any tangent vector fields $X_1, \dots, X_p$ to $M$ we obtain:
\begin{equation*}
\begin{split}
(\phi&\circ\pi)^2\delta_2(\pi^*(\alpha))(X_1,\dots, X_p)=-\sum_{j=1}^{m}(\tilde{f_j}\lrcorner\nabla_{\tilde{f_j}}(\pi^*\alpha))(X_1,\dots, X_p)\\
=&-\sum_{j=1}^{m}\tilde{f_j}(\alpha(\pi_*\tilde{f_j},\pi_*X_1, \dots, \pi_*X_p)\circ\pi)+\sum_{j=1}^{m}\alpha(\pi_*(\nabla_{\tilde{f_j}}\tilde{f_j}), \pi_*X_1,\dots,\pi_*X_p)\circ\pi\\
&+\sum_{j=1}^{m}[
\alpha(\pi_*\tilde{f_j}, \pi_*(\nabla_{\tilde{f_j}}X_1),\dots,\pi_*X_p)+\cdots +
\alpha(\pi_*\tilde{f_j}, \pi_*X_1,\dots,\pi_*(\nabla_{\tilde{f_j}}X_p))]\circ\pi\\
=&\sum_{j=1}^{m}\alpha(\pi_*(\nabla_{\tilde{f_j}}\tilde{f_j}), \pi_*X_1,\dots,\pi_*X_p)\circ\pi\\
=&\sum_{j=1}^{m}\alpha(\pi_*(\widetilde{\nabla^{g_F}_{f_j}f_j}-\frac{g(\tilde{f_j},\tilde{f_j})}{\phi\circ\pi}\mathrm{grad}(\phi\circ\pi)), \pi_*X_1,\dots,\pi_*X_p)\circ\pi\\
=&-m(\phi\circ\pi)(\mathrm{grad}(\phi\circ\pi)\lrcorner\pi^*(\alpha))(X_1,\dots,X_p),
\end{split}
\end{equation*}
where we took into account that $\pi_*(\tilde{f_j})=0$ (as $\tilde{f_j}$ are tangent to the fiber $F$ and $\pi$ is the projection on $B$). This proves \eqref{del2pi}.

For \eqref{del1pi} we compute:
\begin{equation*}
\begin{split}
\delta_1&(\pi^*(\alpha))(X_1,\dots,X_p)=-\sum_{i=1}^{n}(\tilde{e_i}\lrcorner\nabla_{\tilde{e_i}}(\pi^*\alpha))(X_1,\dots, X_p)\\
=&-\sum_{i=1}^{n}\tilde{e_i}(\alpha(\pi_*\tilde{e_i},\pi_*X_1, \dots, \pi_*X_p)\circ\pi)\\
&+\sum_{i=1}^{n}\alpha(\pi_*(\nabla_{\tilde{e_i}}\tilde{e_i}), \pi_*X_1,\dots,\pi_*X_p)\circ\pi\\
&+\sum_{i=1}^{n}[
\alpha(\pi_*\tilde{e_i}, \pi_*(\nabla_{\tilde{e_i}}X_1),\dots,\pi_*X_p)+\cdots +
\alpha(\pi_*\tilde{e_i}, \pi_*X_1,\dots,\pi_*(\nabla_{\tilde{e_i}}X_p))]\circ\pi
\end{split}
\end{equation*}

We may suppose without loss of generality that $X_i$ are lifts of vector fields $Y_i$ on $B$: $X_i=\tilde{Y_i}$, for $i=1,\dots,p$ (since each term in the above sum vanishes if there is at least some vector field $X_i$ tangent to $F$, for which $\pi_*(X_i)$). Under this assumption, we further obtain:
\begin{equation*}
\begin{split}
\delta_1&(\pi^*(\alpha))(X_1,\dots,X_p)=\\
=&-\sum_{i=1}^{n}(e_i(\alpha(e_i,Y_1, \dots, Y_p)))\circ\pi+\sum_{i=1}^{n}\alpha(\nabla^{g_B}_{e_i}e_i, Y_1,\dots,Y_p)\circ\pi\\
&+\sum_{i=1}^{n}[
\alpha(e_i, \nabla^{g_B}_{e_i}Y_1,\dots,Y_p)\circ\pi+\cdots +
\alpha(e_i, Y_1,\dots,\nabla^{g_B}_{e_i}Y_p))\circ\pi]\\
=&(-\sum_{i=1}^{n}e_i\lrcorner\nabla^{g_B}_{e_i}\alpha)(Y_1,\dots,Y_p)\circ\pi=\pi^*(\delta^{g_B}(\alpha))(X_1,\dots,X_p),
\end{split}
\end{equation*}
thus proving \eqref{del1pi}.

Let now $\beta\in\Omega^{p+1}(F)$. We obtain  \eqref{del1sig} as follows:
\begin{equation*}
\begin{split}
\delta_1&(\sigma^*(\beta))(X_1,\dots,X_p)=-\sum_{i=1}^{n}(\tilde{e_i}\lrcorner\nabla_{\tilde{e_i}}(\sigma^*\beta))(X_1,\dots, X_p)\\
=&-\sum_{i=1}^{n}\tilde{e_i}(\beta(\sigma_*\tilde{e_i},\sigma_*X_1, \dots, \sigma_*X_p)\circ\sigma)+\sum_{i=1}^{n}\beta(\sigma_*(\nabla_{\tilde{e_i}}\tilde{e_i}), \sigma_*X_1,\dots,\sigma_*X_p)\\
&+\sum_{i=1}^{n}[
\beta(\sigma_*\tilde{e_i}, \sigma_*(\nabla_{\tilde{e_i}}X_1),\dots,\sigma_*X_p)+\cdots +
\beta(\sigma_*\tilde{e_i}, \sigma_*X_1,\dots,\sigma_*(\nabla_{\tilde{e_i}}X_p))]=0,
\end{split}
\end{equation*}
since $\sigma_*\tilde{e_i}=0$, because $\tilde{e_i}$ is the lift of a vector field on $B$ and also $\sigma_*(\nabla_{\tilde{e_i}}\tilde{e_i})=\sigma_*(\widetilde{\nabla^{g_B}_{e_i}e_i})=0$.

The commutation rule \eqref{del2sig} is shown as follows:
\begin{equation*}
\begin{split}
(\phi&\circ\pi)^2\delta_2(\sigma^*(\beta))(X_1,\dots, X_p)=-\sum_{j=1}^{m}(\tilde{f_j}\lrcorner\nabla_{\tilde{f_j}}(\sigma^*\beta))(X_1,\dots, X_p)\\
=&-\sum_{j=1}^{m}\tilde{f_j}(\beta(\sigma_*\tilde{f_j},\sigma_*X_1, \dots, \sigma_*X_p)\circ\sigma)\\
&+\sum_{j=1}^{m}\beta(\sigma_*(\nabla_{\tilde{f_j}}\tilde{f_j}), \sigma_*X_1,\dots,\sigma_*X_p)\circ\sigma\\
&+\sum_{j=1}^{m}[
\beta(\sigma_*\tilde{f_j}, \sigma_*(\nabla_{\tilde{f_j}}X_1),\dots,\sigma_*X_p)+\cdots +
\beta(\sigma_*\tilde{f_j}, \sigma_*X_1,\dots,\sigma_*(\nabla_{\tilde{f_j}}X_p))]\circ\sigma\\
=&-\sum_{j=1}^{m}f_j(\beta(f_j,\sigma_*X_1, \dots, \sigma_*X_p))\circ\sigma\\
&+\sum_{j=1}^{m}\beta(\sigma_*(\widetilde{\nabla^{g_F}_{f_j}f_j}-\frac{g(\tilde{f_j},\tilde{f_j})}{\phi\circ\pi}\mathrm{grad}(\phi\circ\pi)), \sigma_*X_1,\dots,\sigma_*X_p)\circ\sigma\\
&+\sum_{j=1}^{m}[
\beta(f_j, \sigma_*(\nabla_{\tilde{f_j}}X_1),\dots,\sigma_*X_p)+\cdots +
\beta(f_j, \sigma_*X_1,\dots,\sigma_*(\nabla_{\tilde{f_j}}X_p))]\circ\sigma,
\end{split}
\end{equation*}
where we may again assume, without loss of generality, that $X_i$ are lifts of vector fields $Z_i$ on F: $X_i=\tilde{Z_i}$ for $i=1,\dots,p$. For a tangent vector field $Y$ to $B$, each of the above terms vanishes, since $\sigma_*(Y)=0$. We thus get:

\begin{equation*}
\begin{split}
(\phi&\circ\pi)^2\delta_2(\sigma^*(\beta))(X_1,\dots, X_p)=\\
=&-\sum_{j=1}^{m}f_j(\beta(f_j,Z_1, \dots, Z_p))\circ\sigma+\sum_{j=1}^{m}\beta(\nabla^{g_F}_{f_j}f_j, Z_1,\dots,Z_p)\circ\sigma\\
&+\sum_{j=1}^{m}[\beta(f_j, \nabla^{g_F}_{f_j}Z_1,\dots,\sigma_*X_p)+\cdots +
\beta(f_j, Z_1,\dots,\nabla^{g_F}_{f_j}Z_p)]\circ\sigma\\
=&\sigma^*(\delta^{g_F}(\beta))(X_1,\dots,X_p).
\end{split}
\end{equation*} $\qed$ 
\end{proof}

\begin{Theorem}\label{wp}
Let $(B^n,g_B)$ and $(F^m,g_F)$ be two compact Riemannian manifolds with formal metrics. Then the warped product metric $g=\pi^*(g_B)+(\phi\circ\pi)^2\sigma^*(g_F)$ on $B\times{_\phi} F$ is formal if and only if the warping function $\phi$ is constant.
\end{Theorem}

\begin{proof}
Let $\beta\in\Omega^{p}(F)$ be a $g_F$-harmonic form on $F$ (as $b_m(F)=1$, there exists at least a harmonic $m$-form on $F$).  From \eqref{del1sig} and \eqref{del2sig}, it follows that $\sigma^*\beta$ is a $g$-harmonic form on the warped product $B\times_\phi F$. If we assume the warped metric $g$ to be formal, it follows in particular that the length of $\sigma^*\beta$ is constant. As $g_F$ is also assumed to be formal, the length of $\beta$ is constant as well. On the other hand, the following relation holds:
\begin{equation}\label{len}
g(\sigma^*\beta, \sigma^*\beta)=(\phi\circ\pi)^{2p} g_F(\beta,\beta)\circ\sigma,
\end{equation}
showing that the function $\phi$ must be constant. 

Conversely, if $\phi$ is constant, then the warped product reduces to the Riemannian product between the Riemannian manifolds $(B,g_B)$ and $(F,\phi^2g_F)$, which is geometrically formal by Proposition \ref{prod}. 

$\qed$ 
\end{proof}

\begin{Remark}
From the above proof we see that Theorem \ref{wp} holds more generally for metrics having all harmonic forms of constant length.
\end{Remark}

\section{Geometric formality of Vaisman metrics}

A Vaisman manifold is a particular type of locally conformal K\"ahler (LCK) manifold. It is defined as a Hermitian manifold $(M,J,g)$, of real dimension $n=2m\geq 4$,  whose fundamental $2$-form $\omega$ satisfies the conditions:
$$d\omega=\theta\wedge\omega, \quad \nabla\theta=0.$$
Here $\theta$ is a (closed) $1$-form, called the Lee form, and $\nabla$ is the Levi-Civita connection of the LCK metric $g$ (we always consider $\theta\neq 0$, to not include the K\"ahler manifolds among the Vaisman ones).

Compact Vaisman manifolds are closely related to Sasakian ones, as the following structure theorem shows:

\begin{Theorem}\cite{ov_str}
 Compact Vaisman  manifolds are mapping tori over $S^1$. More precisely: 
 the universal cover $\tilde M$ is a metric cone $N \times \mathbb{R}^{>0}$, with $N$ compact Sasakian manifold
and the deck group is isomorphic with $\mathbb{Z}$, generated by
$(x, t)\mapsto (\lambda(x), t+q)$ for some $\lambda\in \text{Aut}(N)$, $q\in \mathbb{R}^{>0}$.
\end{Theorem}

This puts compact Vaisman manifolds into the framework of warped products and motivates their consideration here.

Vaisman manifolds are abundant. Every Hopf manifold (quotient of $\mathbb{C^N}\setminus\{0\}$ by the cyclic group generated by a semi-simple operator with subunitary eigenvalues) is such, and all its compact complex submanifolds (see \cite[Proposition 6.5]{ve}). Besides, the complete list of Vaisman compact surfaces is given in \cite{be}.

Being parallel and Killing (see \cite{do}), the Lee field $\theta^{\sharp}$ is real holomorphic and, together with $J\theta^{\sharp}$ generates a one-dimensional complex, totally geodesic, Riemannian foliation $\mathcal{F}$. Note that $\mathcal{F}$ is transversally K\"ahler. 

In the sequel, the terms {\em basic (foliate)} and {\em horizontal} refer to $\mathcal{F}$. Moreover, we shall use the basic versions of the standard operators acting on $\Omega^*_B(M)$, the space of basic forms: $\Delta_B$ is the basic Laplace operator, $L_B$ is the exterior multiplication with the transversal K{\"a}hler form and $\Lambda_B$ its adjoint with respect to the transversal metric.

The main result of this section puts severe restrictions on formal Vaisman metrics:

\begin{Theorem}\label{gfvais}
Let $(M^{2m},g,J)$ be a compact Vaisman manifold. The metric $g$ is geometrically formal if and only if $b_p(M)=0$ for $2\leq p\leq 2m-2$ and $b_1(M)=b_{2m-1}(M)=1$.
\end{Theorem}

\begin{proof}
Let $\gamma\in\Omega^p(M)$ be a harmonic form on $M$ for some $p$, $1\leq p\leq m-1$. By \cite[Theorem~4.1]{vaism}, $\gamma$ has the following form:
\begin{equation}\label{decomp}
\gamma=\alpha+\theta\wedge \beta,
 \end{equation}
with $\alpha$ and $\beta$ basic, transversally harmonic and transversally primitive.

Since $\alpha$ is basic, $J\alpha$ is also a basic $p$-form that is transversally harmonic and transversally primitive:
\[\Delta_B(J\alpha)=0, \quad \Lambda_B(J\alpha)=0,\]
because $\Delta_B$ and $\Lambda_B$ both commute with the transversal complex structure $J$ (as the foliation is transversally K\"ahler). 
Again from \cite[Theorem~4.1]{vaism}, by taking $\beta=0$, it follows that $J\alpha$ is a harmonic form on $M$: $\Delta(J\alpha)=0$.

The assumption that $g$ is geometrically formal implies that $\alpha \wedge J\alpha$ is harmonic on $M$, so that in particular it is coclosed: $\delta(\alpha \wedge J\alpha)=0$. According to \cite{vaism}, this implies that $\alpha \wedge J\alpha$ is transversally primitive: $\Lambda_B(\alpha \wedge J\alpha)=0$. 

On the other hand, it is proven in \cite[Proposition~2.2]{grosnagy} that 
for primitive forms $\eta,\mu\in\Lambda^p V$, where $(V,g,J)$ is any Hermitian vector space, the following algebraic relation holds:
\begin{equation}\label{algrel}
(\Lambda)^p(\eta\wedge\mu)=(-1)^{\frac{p(p-1)}{2}}p!\<\eta,J\mu>,
\end{equation}
where $J$ is the extension of the complex structure to $\Lambda^*V$ defined by:
\[(J\eta)(v_1,\dots,v_p):=\eta(Jv_1,\dots,Jv_p), \quad \text{for all } \eta\in\Lambda^pV, v_1,\dots,v_p\in V.\]
We apply the above formula to the transversal K\"ahler geometry and obtain  that $\alpha$ vanishes everywhere:
\[0=(\Lambda_B)^p(\alpha \wedge J\alpha)=(-1)^{\frac{p(p+1)}{2}}p!\<\alpha,\alpha>.\]

The same argument as above applied to $\beta\in\Omega^{p-1}_B(M)$ shows that $\beta$ is identically zero if $p\geq2$. Thus, $\gamma=0$ for $2\leq p\leq m-1$, which proves that:
\[b_2(M)=\cdots=b_{m-1}(M)=0.\]

If  $p=1$, then $\beta$ is a basic function, which is transversally harmonic, so that $\beta$ is a constant. Thus $\gamma$ is a multiple of $\theta$, showing that the space of harmonic $1$-forms on $M$ is $1$-dimensional: $b_1(M)=1$. 

It remains to show that the Betti number in the middle dimension, $b_m(M)$, also vanishes. This follows from Proposition~\ref{bettzero} applied to $p=1$, $k=m$ and $l=2$.

The converse is clear, since the space of harmonic forms with respect to the Vaisman metric $g$ is spanned by $\{1,\theta,*\theta,dvol_g\}$ and thus the only product of harmonic forms which is not trivial is $\theta\wedge *\theta=g(\theta,\theta)dvol_g$, which is harmonic because $\theta$ has constant length, being a parallel $1$-form.

$\qed$ 
\end{proof}

\begin{Remark}
Theorem~\ref{gfvais} may be considered as an analogue of the following result on the geometric formality of Sasakian manifolds:
\begin{Theorem}\cite[Theorem~2.1]{grosnagy} \label{gfsas}
Let $(M^{2n+1},g)$ be a compact Sasakian manifold. If the metric $g$ is geometrically formal, then $b_p(M)=0$ for $1\leq p\leq 2n$, \emph{i.e.} $M$ is a real cohomology sphere.
\end{Theorem}
\end{Remark}

\noindent{\bf Acknowledgement.} We thank D. Kotschick and P.-A. Nagy for very useful comments on the first draft of this note.

\end{document}